\documentclass[12pt,reqno]{amsart}
\usepackage{amssymb}
\usepackage{amsthm}
\usepackage{amsmath}
\usepackage{a4wide}
\usepackage{latexsym}
\usepackage[v2,tips]{xy}
\usepackage{mathrsfs}
\usepackage[T1]{fontenc}
\usepackage[latin1]{inputenc}
\usepackage{enumitem}
\usepackage{graphicx}
\usepackage{hyperref}
\usepackage[normalem]{ulem}
\usepackage{dutchcal}
\usepackage{tikz-cd} 
\usepackage{cleveref}
\usepackage{tabularx}

\newcommand{\N}{\mathbb{N}}

\newcommand{\R}{\mathbb{R}}
\newcommand{\C}{\mathbb{C}}
\newcommand{\K}{\mathbb{K}}

\newtheorem{theorem}{Theorem}[section]
\newtheorem{lemma}[theorem]{Lemma}
\newtheorem{proposition}[theorem]{Proposition}
\newtheorem{corollary}[theorem]{Corollary}
\theoremstyle{definition}
\newtheorem{definition}[theorem]{Definition}

\newtheorem{remark}[theorem]{Remark}

\newtheorem*{organization}{Organization and overview of content}
\newtheorem*{acknowledgements}{Acknowledgements}

\numberwithin{equation}{section}

\renewcommand{\ge}{\geqslant}

\renewcommand{\le}{\leqslant}
\renewcommand{\leq}{\leqslant}

\newcommand{\gl}{\ensuremath{\Gamma}\!\operatorname{L}_1}
\newcommand{\SC}{\operatorname{SC}}

\title[Strictly singular operators on Baernstein and Schreier spaces]{Strictly singular operators on\\ the Baernstein and Schreier spaces}
\author[N.J.~Laustsen]{Niels Jakob Laustsen}
\address{(N.J.~Laustsen) School of Mathematical Sciences, Fylde
  College, Lancaster University, Lancaster LA1 4YF, United Kingdom}
\email{n.laustsen@lancaster.ac.uk}

\author[J.~Smith]{James Smith}
\address{(J.~Smith) School of Mathematical Sciences, Fylde
  College, Lancaster University, Lancaster LA1 4YF, United Kingdom}
\email{j.smith43@lancaster.ac.uk}

\subjclass{46H10,
47L10,
(primary); 
46B03,
46B45,   	
47L20 (secondary)}
\keywords{Banach space, Baernstein space, Schreier space, bounded operator, strictly singular operator, compact operator, subsymmetric basic sequence, uncomplemented block subspace}

\begin{document}

\begin{abstract}
Every composition of two strictly singular operators is compact on the Baernstein space~$B_p$ for $1 < p < \infty$ and on the $p$-convexified Schreier space~$S_{p}$ for $1 \leq p < \infty$. Furthermore, every subsymmetric basic sequence in~$B_p$ (respectively,~$S_p$) is equivalent to the unit vector basis for~$\ell_p$ (respectively,~$c_0$), and the Banach spaces~$B_p$ and~$S_p$ contain block basic sequences whose closed span is not complemented. 
\end{abstract}

\maketitle

\section{Introduction}
\noindent
This note is a continuation of the line of research we initiated in~\cite{LS}, investigating the Baernstein and $p$-convexified Schreier spaces, particularly their lattices of closed ideals of bounded operators. Two of the main conclusions of~\cite{LS} are that on each of these spaces, there are~$2^{\mathfrak{c}}$ closed ideals lying between the ideals of compact and strictly singular operators, as well as $2^{\mathfrak{c}}$ closed ideals that contain infinite-rank projections. 

These results are counterparts of, and were motivated by, the seminal work~\cite{JS} of John\-son and Schechtman concerning closed ideals of bounded operators on the Lebesgue spaces $L_p[0,1]$ for $p\in(1,2)\cup(2,\infty)$.  Johnson and Schechtman~\cite{JS2} went on to show that among all these ideals, only the ideal of compact operators possesses any kind of approximate identity. (Monotonicity of the basis and reflexivity ensure that the sequence of basis projections is a contractive, two-sided approximate identity for the ideal of compact operators on~$L_p[0,1]$.)

This raises the question: Which closed ideals of bounded operators on the Baernstein and Schreier spaces possess approximate identities? We shall answer it for the ``small'' ideals (those that are contained in the ideal of strictly singular operators), obtaining the same conclusion as Johnson and Schechtman; that is, only the ideal of compact operators does. The reason is that the composition of two strictly singular operators on any Baernstein or $p$-convexified Schreier space is compact. 

Spaces with this property are well known in the literature. Probably the oldest and most famous instances are the continuous functions~$C(K)$ on a compact Hausdorff space~$K$ and the Lebesgue spaces~$L_1(\mu)$ for a $\sigma$-finite measure~$\mu$. 
This follows by combining two classical results: 
\begin{itemize}
    \item $C(K)$ and~$L_1(\mu)$ have the Dunford--Pettis property (see for example \cite[Theorem~5.4.5]{AK} or \cite[Theorems~VI.7.4 and~VI.8.12]{DS}); this implies that the composition of two weakly compact operators is compact.
    \item A bounded operator on~$C(K)$ or~$L_1(\mu)$ is strictly singular if and only if it is weakly compact (this is due to Pe\l{}czy\'{n}ski~\cite{P1,P2}).
\end{itemize}
Milman~\cite{milman} complemented the $L_1$-case by showing that the composition of two strictly singular operators on~$L_p[0,1]$ is also compact when $1<p<\infty$. 

Much closer to our investigation, Causey and Pel\-czar-Bar\-wacz \cite[Theorem~7.6, (3)--(4)]{CPB} have recently studied compositions of strictly singular operators on the higher-order Schreier spaces~$X[\mathcal{S}_\xi]$ for every countable ordinal $\xi>0$, proving that there is a (necessarily unique) natural number~$k$, dependent only on the ordinal~$\xi$, such that:
\begin{itemize}
\item every composition of $k+1$ strictly singular operators on~$X[\mathcal{S}_\xi]$ is compact;
\item there are~$k$ strictly singular operators on~$X[\mathcal{S}_\xi]$ whose composition is not compact.
\end{itemize}
We note that $k=\xi$ when~$\xi$ is finite, so in particular $k=1$ for~$X[\mathcal{S}_1]$, which is the space we denote~$S_1$. As indicated above, our main result extends this conclusion to the $p$-convexified variants of~$S_1$, as well as the Baernstein spaces~$B_p$ for $1<p<\infty$. We now state it formally.
\begin{theorem}\label{StimesSisK}
    Let~$T$ and~$U$ be strictly singular operators on either the Baernstein space~$B_p$ for some $1 < p < \infty$ or the $p$-convexified Schreier space~$S_{p}$ for some $1 \leq p < \infty$. Then the composite operator $TU$ is compact.    
\end{theorem}

\begin{organization} In \Cref{S:2}, we introduce the set-up and main objects that we study, notably giving the precise definitions of the Banach spaces~$S_p$ and~$B_p$, as well as the Gas\-pa\-ris--Leung index, which we then use to prove some results about domination and equivalence of subsequences of normalized block basic sequences of the unit vector bases for~$S_p$ and~$B_p$.

These results lay the foundations for our proof of \Cref{StimesSisK}, which we present in \Cref{S:3}. In fact, \Cref{S:3} contains \emph{two} proofs of it: The first is self-contained, while the other relies on the theory of Schreier spreading basic sequences and equivalence thereof, which is due to Androulakis \textit{et al}~\cite{ADST}. However, it also depends on some of the lemmas that we established in the course of our first proof, so the two proofs are not truly independent. 

We conclude with two short sections. In \Cref{S:4}, we apply one of our results about basic sequences to deduce that the only subsymmetric basic sequences in~$B_p$ (respectively,~$S_p$) are those equivalent to the unit vector basis for~$\ell_p$ (respectively,~$c_0$). Finally, in \Cref{S:5}, we show that~$B_p$ and~$S_p$ contain block basic sequences whose closed span is not complemented. 
\end{organization}

\section{Block basic sequences in the Baernstein and Schreier spaces}\label{S:2}

\noindent Before proving our first results, we state the main definitions and conventions regarding notation and terminology, all of which are similar to~\cite{LS}. In particular, $\K$~denotes the scalar field, either~$\R$ or~$\C$, and we use function notation for sequences, so that~$x(n)$ is the $n^{\text{th}}$ coordinate of a sequence~$x\in\K^\N$. As usual, $\operatorname{supp} x = \{ n\in\N : x(n)\ne 0\}$ for $x\in\K^\N$, and~$c_{00}$ denotes the subspace of finitely supported elements of~$\K^\N$. It is spanned by the ``unit vector basis'' $(e_{n})_{n\in\N}$ given by $e_n(m) = 1$ if $m=n$ and $e_n(m) = 0$ otherwise.  We shall sometimes consider two Banach spaces~$D$ and~$E$ simultaneously, both of which have the unit vector basis as a (Schauder) basis. In such instances, we write $(d_n)_{n\in\N}$ for the unit vector basis of~$D$ and $(e_n)_{n\in\N}$ for the unit vector basis of~$E$ to keep track of the ambient spaces. 

By an \emph{operator,} we understand a bounded linear map between Banach spaces. We write $\mathscr{B}(X,Y)$ for the Banach space of operators from a Banach space~$X$ into a Banach space~$Y$, abbreviated $\mathscr{B}(X)$ when $X=Y$; $I_X$~denotes the identity operator on~$X$.  An operator $T\in\mathscr{B}(X,Y)$ is \emph{strictly singular} if no restriction of~$T$ to an infinite-dimensional subspace of~$X$ is an isomorphic embedding. 

Our focus is on two classes of Banach spaces: the $p$-convexified Schreier spaces and the Baern\-stein spaces. The cornerstone of their definitions is the notion of a
\emph{Schreier set,} that is, a finite subset~$F$ of the natural numbers $\N=\{1,2,3,\ldots\}$ such that either $F=\emptyset$ or $\lvert F\rvert\leq \min F$, where $\lvert F\rvert$ denotes the cardinality of~$F$. In line with standard practice, $\mathcal{S}_{1}$ denotes the family of Schreier sets. A frequently used property of~$\mathcal{S}_{1}$ is that it is \emph{spreading} in the sense that if $F=\{f_1<\cdots<f_n\}$ belongs to~$\mathcal{S}_1$ and $G=\{g_1<\cdots<g_n\}$ is another subset of~$\N$ which satisfies $f_j\le g_j$ for each $1\le j\le n$, then also $G\in\mathcal{S}_1$. 

For $1\le p<\infty$ and $F\in\mathcal{S}_1$, we can define a seminorm~$\mu_p(\,\cdot\,,F)$ on~$\K^\N$ by
\[ \mu_{p}(x,F) =\begin{cases} 0\ &\text{if}\ F=\emptyset,\\ \displaystyle{\biggl(\sum_{n \in F}\lvert x(n)\rvert^{p}\biggr)^{1/p}}\ &\text{otherwise} \end{cases}\qquad (x\in\K^\N). \] 
Set $\lVert x\rVert_{S_{p}} = \sup\{ \mu_{p}(x,F) : F \in \mathcal{S}_{1}\}\in[0,\infty]$ for $x\in\K^\N$ and \mbox{$Z_p =\{ x\in\K^\N : \lVert x \rVert_{S_{p}} < \infty\}$.} Standard arguments show that~$Z_p$ is a subspace of~$\K^\N$ and~$\lVert\,\cdot\,\rVert_{S_p}$ a complete norm on it. Furthermore, by \cite[Propositions~3.5 and~3.10 and Corollary~3.12]{BL}, the unit vector basis $(e_{n})_{n\in\N}$ is a $1$\nobreakdash-un\-con\-di\-tional, shrinking, normalized basic sequence in~$Z_p$ and hence a basis for its closed linear span, which we denote~$S_p$ and call the \emph{$p$-convexified Schreier space}. We remark that~$S_p$ is a proper subspace of~$Z_p$; in fact, the latter is non-separable, as observed in \cite[Corollary~5.6]{BL}. 

To analogously define the Baernstein spaces, we require the notion of a \emph{Schreier chain}, which is a non-empty, finite collection~$\mathcal{C}$ of non-empty, consecutive Schreier sets; that is, 
\[ \mathcal{C} = \{F_{1},\ldots,F_{m} :  m\in\N,\, F_1,\ldots,F_m\in\mathcal{S}_1\setminus\{\emptyset\},\, F_1<F_2<\cdots<F_m\}, \] 
using the standard notational convention that $F_j<F_{j+1}$ means $\max F_{j} < \min F_{j+1}$ for (finite, non-empty) subsets~$F_j $ and~$F_{j+1}$ of~$\N$. 

Given $1<p<\infty$ and a Schreier chain~$\mathcal{C}$, we can define a seminorm $\beta_p(\,\cdot\,,\mathcal{C})$ on~$\K^\N$ by 
\[ \beta_{p}(x,\mathcal{C}) = \biggl(\sum_{F\in \mathcal{C}}\Bigl(\sum_{n\in F}\lvert x(n)\rvert\Bigr)^{p}\biggr)^{1/p} = \biggl(\sum_{F\in \mathcal{C}}\mu_1(x,F)^p\biggr)^{1/p}\qquad (x\in\K^\N). \]
Set $\lVert x\rVert_{B_{p}} = \sup\{\beta_{p}(x,\mathcal{C}) : \mathcal{C}\in\SC\}\in[0,\infty]$ for $x\in\K^\N$, where~$\SC$ denotes the collection of all Schreier chains. 
Then $B_p =\{ x\in\K^\N : \lVert x \rVert_{B_{p}} < \infty\}$ is a subspace of~$\K^\N$ and~$\lVert\,\cdot\,\rVert_{B_p}$ a complete norm on it. In contrast to the Schreier spaces, $c_{00}$ is dense in~$B_p$, which is the $p^{\text{th}}$ \emph{Baern\-stein space}. It is reflexive, and the unit vector basis $(e_{n})_{n\in\N}$ is a $1$\nobreakdash-un\-con\-di\-tional, normalized basis for it. Baernstein~\cite{AB} originally defined the space~$B_2$, while Seifert~\cite{S} ob\-served that Baernstein's definition generalizes to arbitrary~$p>1$.

As in~\cite{LS}, we shall often consider the Baernstein and Schreier spaces simultaneously, using the letter~$E$ to denote either~$B_p$ for some $1<p<\infty$ or~$S_p$ for some $1\le p<\infty$. Then, for a subset~$N$ of~$\N$, we write~$E_N$ for the closed subspace of~$E$ spanned by $\{e_n : n\in N\}$. 
This subspace is $1$-complemented in~$E$ by $1$-unconditionality of the unit vector basis. We write $(e_n^*)_{n\in\N}$ for the sequence of coordinate functionals associated with the unit vector basis. They form a $1$-unconditional basis for the dual space~$E^*$ because the unit vector basis for~$E$ is shrinking.  

\smallskip 

Gasparis and Leung~\cite{GL} introduced a numerical index~$\gl(M,N)$ for every pair $(M,N)$ of infinite subsets of~$\N$ to help analyze the closed subspaces spanned by subsequences of the unit vector basis for the Schreier space~$S_1$ and its higher-order counterparts. We showed in~\cite[Section~4]{LS} that this index has similar properties for the Baernstein and Schreier spaces~$B_p$ and~$S_p$ for $p>1$. In the context of the present investigation, its most important property is that for infinite subsets~$M$ and~$N$ of~$\N$, the basic sequence $(e_m)_{m\in M}$ dominates the basic sequence $(e_n)_{n\in N}$ in either~$B_p$ or~$S_p$ if and only if
$\gl(M,N)<\infty$. 

In order to state concisely \cite[Definition~3.3]{GL}, in which the index~$\gl(M,N)$ is defined, we require two pieces of notation. First, for an infinite subset $M = \{m_1<m_2<\cdots\}$ of~$\mathbb{N}$, set 
\[ M(J) = \{m_{j} : j \in J\}\qquad (J \subseteq \mathbb{N}). \] 
Second, the \emph{Schreier covering number} of a finite subset~$A$ of~$\N$ is given by $\tau_1(\emptyset) = 0$ and 
\[ \tau_1(A) = \min\biggl\{ m\in \mathbb{N} : A \subseteq \bigcup_{j=1}^{m}F_{j},\ \text{where}\ F_1,\ldots,F_m\in\mathcal{S}_1\ \text{and}\ F_1<F_2<\cdots<F_m\biggr\} \]
for $A\ne\emptyset$. 
Now the \emph{Gas\-pa\-ris--Leung index} of an infinite subset $M$ of~$\N$ with respect to another infinite subset~$N$ of~$\N$ is 
\begin{equation}\label{Eq:GLindex}
\gl(M,N) = \sup\bigl\{\tau_1(M(J)) : J\subset\mathbb{N}\ \text{finite,}\ N(J) \in \mathcal{S}_{1}\bigr\}\in\N\cup\{\infty\}.
\end{equation}
(We remark that Gasparis and Leung denoted this index~$d_1(M,N)$. We have chosen the more distinctive symbol~$\gl(M,N)$ in their honour, noting that~$\Gamma$ is the first letter of the Greek spelling of ``Gasparis''.)

\begin{remark}\label{R:tau} For $n\in\N$ and subsets~$M$ and~$J$ of~$\N$, where~$M$ is infinite and~$J$ is finite and non-empty, we have  $\tau_1(M(J))\le n$ if and only if we can find a natural number $m\le n$ and Schreier sets $F_1<\cdots<F_m$ such that $M(J)=\bigcup_{j=1}^{m}F_{j}$. Writing $F_j = M(J_j)$, we see that $\tau_1(M(J))\le n$ if and only if we can decompose~$J$ as $J = \bigcup_{j=1}^m J_j$ for some natural number $m\le n$, where the sets $J_1<\cdots <J_m$ satisfy $M(J_j)\in\mathcal{S}_1$ for each $1\le j\le m$.  
\end{remark}

We shall now establish some basic properties of the Gasparis--Leung index. We state them in greater generality than is strictly necessary for the applications we have in mind, for two reasons. First, we expect that they will be useful in future studies of the Baernstein and Schreier spaces, and second, we hope that they may inspire others to continue this line of research, investigating the Gasparis--Leung index from a purely com\-bi\-na\-torial/number-theo\-retic point of view.

\begin{lemma}\label{L:GLindexProperties}
  Let $L$, $M = \{m_1<m_2<\cdots\}$ and $N = \{n_1<n_2<\cdots\}$ be infinite subsets of~$\N$. Then:
  \begin{enumerate}[label={\normalfont{(\roman*)}}]  
  \item\label{L:GLindexProperties1} $\gl(L,N)\le \gl(L,M)\cdot\gl(M,N)$.
  \item\label{L:GLindexProperties2} Suppose that~$N$ is a spread of~$M$ in the sense that $n_j\ge m_j$ for every $j\in\N$. Then \[ \gl(N,M)=1\qquad\text{and}\qquad \gl(L,M)\le\gl(L,N). \] 
  \item\label{L:GLindexProperties3} $\gl(N,N\setminus F)\le \tau_1(\{n_j : 1\le j\le \lvert F\rvert\}) + 1$ for every finite subset~$F$ of~$N$.
  \item\label{L:GLindexProperties4} Suppose that $m_j\le n_{j+1}$ for each $j\in\N$. Then $\gl(N,M)\le 2$.
\end{enumerate}    
\end{lemma}

\begin{proof}
    \ref{L:GLindexProperties1}. Take a non-empty, finite subset~$J$ of~$\N$ for which $N(J)\in\mathcal{S}_1$. Then we have $\tau_1(M(J))\le\gl(M,N)$, so \Cref{R:tau} implies that we can write $J = \bigcup_{j=1}^m J_j$ for some natural number $m\le\gl(M,N)$, where the sets $J_1<\cdots <J_m$ satisfy $M(J_j)\in\mathcal{S}_1$ for each $1\le j\le m$. Applying \Cref{R:tau} again, for each $1\le j\le m$, we can find a natural number $r_j\le\gl(L,M)$ and sets $J_{j,1}<J_{j,2}<\cdots <J_{j,r_j}$ such that $J_j = \bigcup_{k=1}^{r_j} J_{j,k}$ and $L(J_{j,k})\in\mathcal{S}_1$ for each $1\le k\le r_j$. It follows that $J = \bigcup_{j=1}^m\bigcup_{k=1}^{r_j} J_{j,k}$, where 
    \[ J_{1,1}<J_{1,2}<\cdots <J_{1,r_1} < J_{2,1} < J_{2,2} < \cdots < J_{2,r_2} < \cdots < J_{m,1} < \cdots < J_{m,r_m}, \]
    so appealing to \Cref{R:tau} once more, we conclude that
    \[ \tau_1(L(J))\le\sum_{j=1}^m r_j\le m\cdot\max_{1\le j\le m}r_j\le\gl(M,N)\cdot\gl(L,M). \]

    \ref{L:GLindexProperties2}. Suppose that~$N$ is a spread of~$M$. Then~$N(J)$ is a spread of~$M(J)$ for every \mbox{$J\subseteq\N$}, so in particular $N(J)\in\mathcal{S}_1$ whenever $M(J)\in\mathcal{S}_1$. This proves that $\gl(N,M)=1$, and also that $\tau_1(L(J))\le \gl(L,N)$ for every $J\subset\N$ such that $M(J)\in\mathcal{S}_1$. The second inequality follows.

    \ref{L:GLindexProperties3}. Set $k=\lvert F\rvert$, $G=\{n_j : 1\le j\le k\}$ and $N' = N\setminus G= \{n_{k+j} : j\in\N\}$. Then~$N'$ is a spread of $N\setminus F$, so by~\ref{L:GLindexProperties2}, it suffices to show that $\gl(N,N')\le \tau_1(G)+1$; that is, $\tau_1(N(J))\le \tau_1(G)+1$ for every non-empty, finite subset $J=\{j_1<\cdots<j_m\}$ of~$\N$ such that $N'(J)\in\mathcal{S}_1$. Set $K = \{ j_i : 1\le i\le\min\{k,m\}\}$, which is a spread of the inter\-val $[1,\min\{k,m\}]\cap\N$. Con\-se\-quent\-ly~$N(K)$ is a spread of $\{n_i : 1\le i\le\min\{k,m\}\}$, which is a sub\-set of~$G$, and hence $\tau_1(N(K))\le \tau_1(G)$. This completes the proof if $J=K$. Other\-wise $k<m$; then $\max K = j_k$, so $\min (J\setminus K) = j_{k+1}$, and therefore
    \begin{alignat*}{2} 
    \min N(J\setminus K) &= n_{j_{k+1}}\ge n_{k+j_1} = \min N'(J)\ge \lvert N'(J)\rvert\qquad &&\text{because}\ N'(J)\in\mathcal{S}_1\\ 
    &=\lvert J\rvert > \lvert J\setminus K\rvert = \lvert N(J\setminus K)\rvert.  \end{alignat*}
    This proves that $N(J\setminus K)\in\mathcal{S}_1$. Since $K<J\setminus K$, we conclude that 
    \[ \tau_1(N(J))\le \tau_1(N(K)) + 1\le \tau_1(G)+1, \]
    as desired.

    \ref{L:GLindexProperties4}. We have
    \[ \gl(N,M) \le\gl(N,N\setminus\{n_1\})\cdot \gl(N\setminus\{n_1\},M)\le \bigl(\tau_1(\{n_1\})+1\bigr)\cdot 1 = 2, \]
    where the first inequality follows from~\ref{L:GLindexProperties1} and the second from~\ref{L:GLindexProperties3} and~\ref{L:GLindexProperties2} because the hypo\-thesis implies that $N\setminus\{n_1\}$ is a spread of~$M$.     
\end{proof}

\begin{corollary}\label{C:intertwiningSeqs}
  Let $(m_j)_{j\in\N}$ and $(n_j)_{j\in\N}$ be increasing sequences of natural numbers, and let $(e_j)_{j\in\N}$ denote the unit vector basis for~$E$, where $E=B_p$ for some $1<p<\infty$ or $E=S_p$ for some $1\le p<\infty$.  
  \begin{enumerate}[label={\normalfont{(\roman*)}}]  
  \item\label{C:intertwiningSeqs1} Suppose that $m_j\le n_j$ for each $j\in\N$. Then the basic sequence $(e_{n_j})_{j\in\N}$ $1$-dominates $(e_{m_j})_{j\in\N}$.
  \item\label{C:intertwiningSeqs2} Suppose that $m_j\le n_{j+1}$ for each $j\in\N$.  Then the basic sequence $(e_{n_j})_{j\in\N}$ $C$-dom\-i\-nates $(e_{m_j})_{j\in\N}$, where $C=2$ for $E=B_p$ and $C=2^{1/p}$ for $E=S_p$. 
  \item\label{C:intertwiningSeqs3} There exists an infinite subset~$J$ of~$\N$ such that the basic sequences $(e_{m_j})_{j\in J}$ and $(e_{n_j})_{j\in J}$ are equivalent. 
  \end{enumerate}
\end{corollary}

\begin{proof} \ref{C:intertwiningSeqs1}--\ref{C:intertwiningSeqs2}. According to \Cref{L:GLindexProperties}\ref{L:GLindexProperties2} and~\ref{L:GLindexProperties4}, the sets $M=\{m_1<m_2<\cdots\}$ and $N=\{n_1<n_2<\cdots\}$ satisfy $\gl(N,M)=1$ in case~\ref{C:intertwiningSeqs1} and $\gl(N,M)\le 2$ in case~\ref{C:intertwiningSeqs2}. Now the conclusions follow from \cite[Lemma~4.2]{LS}. 

\ref{C:intertwiningSeqs3}. By recursion, we can choose integers $1=j_1<j_2<j_3<\cdots$ such that $m_{j_k}\le n_{j_{k+1}}$ and $n_{j_k}\le m_{j_{k+1}}$ for each $k\in\N$. Then~\ref{C:intertwiningSeqs2} implies that the subsequences $(e_{m_{j_k}})_{k\in\N}$ and $(e_{n_{j_k}})_{k\in\N}$ $C$-dominate each other, so the subset $J = \{j_1<j_2<\cdots\}$ of~$\N$ has the required property.
\end{proof}

The distinction between ``flat'' and ``spiky'' vectors plays a key role in the study of the Baernstein and Schreier spaces. The uniform norm given by $\lVert x\rVert_\infty = \sup_{n\in\N}\lvert x(n)\rvert\in[0,\infty]$ for $x\in\K^\N$ is the tool that allows us to quantify this distinction. (Note that~$B_p$ and~$S_p$ are contained  in~$c_0$ as vector subspaces because the unit vector basis is a normalized basis for them. Hence $\lVert x\rVert_\infty<\infty$ for~$x$ belonging to any of these spaces.)

The importance of certain vectors being ``spiky'' is already evident in the main theorem \cite[Theorem~1.1]{GL} of Gasparis and Leung (who denoted the uniform norm~\mbox{$\lVert\,\cdot\,\rVert_0$}, not~\mbox{$\lVert\,\cdot\,\rVert_\infty$}) and in our counterpart \cite[Theorem~4.1]{LS} of it for the Baernstein and $p$-convexified Schreier spaces. 

On the other hand, ``flat'' block basic sequences appear in \cite[Lemma~3.14]{GL} and \cite[Proposition~2.14]{LS}; we include the details of the latter result for later reference and to provide context for our next lemma, whose final part is a counterpart of it for ``spiky'' block basic sequences. 

\begin{proposition}[{\cite[Proposition~2.14]{LS}}]\label{P:charDsubseq}
    Let $(E,D) = (B_{p},\ell_{p})$ for some $1 < p < \infty$ or $(E,D) = (S_{p},c_{0})$ for some $1 \leq p < \infty$. The following conditions are equivalent for a normalized block basic sequence $(u_{n})_{n\in\N}$  of the unit vector basis for~$E\colon$
    \begin{enumerate}[label={\normalfont{(\alph*)}}]    
    \item\label{P:charDsubseq:a} $\inf_{n\in\N}\lVert u_n\rVert_\infty = 0;$
    \item\label{P:charDsubseq:b} $(u_{n})_{n\in\N}$ admits a subsequence which is $C$-equivalent to the unit vector basis for~$D$, for every constant $C>1;$
    \item\label{P:charDsubseq:c} $(u_{n})_{n\in\N}$ admits a subsequence which is dominated by the unit vector basis for~$D$. 
    \end{enumerate}
\end{proposition}

\begin{proposition}\label{L:nbbs} 
Let $(u_j)_{j\in\N}$ be a normalized block basic sequence of the unit vector basis $(e_j)_{j\in\N}$ for~$E$, where $E=B_p$ for some $1<p<\infty$ or $E=S_p$ for some $1\le p<\infty$. 
\begin{enumerate}[label={\normalfont{(\roman*)}}]  
\item\label{L:nbbs1} Set $m_j = \max\operatorname{supp} u_j$ for $j\in\N$. Then $(e_{m_j})_{j\in\N}$ $C_1$-dom\-i\-nates $(u_j)_{j\in\N}$, where $C_1=3^{1/p}$ if $E=B_p$ and $C_1=1$ if $E=S_p$.
\item\label{L:nbbs2} Suppose that $\delta :=\inf_{j\in\N}\lVert u_j\rVert_\infty>0$, and choose $n_j\in\operatorname{supp} u_j$ such that \mbox{$\lvert\langle u_j, e_{n_j}^*\rangle\rvert\ge\delta$} for each $j\in\N$. Then $(u_j)_{j\in\N}$ $\delta^{-1}$-dominates $(e_{n_j})_{j\in\N}$, and $\overline{\operatorname{span}}\,\{ u_j : j\in\N\}$ is $C_2$\nobreakdash-com\-ple\-mented in~$E$, where $C_2=2\cdot 3^{1/p}/\delta$ if $E=B_p$ and $C_2=2^{1/p}/\delta$ if $E=S_p$.
\item\label{L:nbbs3} The following conditions are equivalent:  
  \begin{enumerate}[label={\normalfont{(\alph*)}}]  
    \item\label{P:pointedbbs1} $\inf_{j\in\N}\lVert u_j\rVert_\infty>0;$
    \item\label{P:pointedbbs2} $(u_j)_{j\in\N}$ is equivalent to a subsequence of~$(e_j)_{j\in\N};$
    \item\label{P:pointedbbs3} $(u_j)_{j\in\N}$ dominates a subsequence of~$(e_j)_{j\in\N}$.
  \end{enumerate}
\end{enumerate}    
\end{proposition}

\begin{proof} Before we embark on the proof, let us record that the sequences $(m_j)_{j\in\N}$ and $(n_j)_{j\in\N}$ in parts~\ref{L:nbbs1} and~\ref{L:nbbs2} are increasing because $(u_j)_{j\in\N}$ is a block basic sequence. 

\ref{L:nbbs1}. Set $M = \{ m_j:j\in\N\}$. Our aim is to prove that for $k\in\N$ and $\alpha_1,\ldots,\alpha_k\in\K$, the elements $x = \sum_{j=1}^k\alpha_j e_{m_j}$ and $y=\sum_{j=1}^k\alpha_j u_{j}$ satisfy the inequality $C_1\lVert x\rVert_E\ge \lVert y\rVert_E$. This is trivial if $\alpha_1=\cdots=\alpha_k=0$, so we may suppose otherwise; thus $x,y\ne 0$. 

 If $E=S_p$, take $F\in\mathcal{S}_1$ such that $\lVert y\rVert_{S_p} = \mu_p(y,F)$. After replacing~$F$ with $F\cap\operatorname{supp} y$, we may suppose that $F\subseteq\operatorname{supp} y$. Set $F_j = F\cap\operatorname{supp} u_j\in\mathcal{S}_1$ for $1\le j\le k$. Then, defining $J = \{ j\in\{1,\ldots,k\} : F_j\ne \emptyset\}$ and  $j_1=\min J$, we have
\begin{alignat*}{2} \min M(J) = m_{j_1} &\ge\min F\qquad &&\text{because}\ \min F\in\operatorname{supp} u_{j_1}\\
&\ge\lvert F\rvert &&\text{because}\ F\in\mathcal{S}_1\\
&\ge\lvert J\rvert &&\text{because}\ (F_j)_{j\in J}\ \text{are disjoint, non-empty subsets of}\ F\\
&= \lvert M(J)\rvert,
\end{alignat*}
so $M(J)\in\mathcal{S}_1$. Consequently,
\begin{align*}
    \lVert x\rVert_{S_p}^p &\ge \mu_p(x,M(J))^p = \sum_{j\in J}\lvert\alpha_j\rvert^p = \sum_{j\in J}\lvert\alpha_j\rvert^p\lVert u_j\rVert_{S_p}^p\\ &\ge \sum_{j\in J}\lvert\alpha_j\rvert^p\mu_p(u_j,F_j)^p
   = \sum_{j\in J}\mu_p(y,F_j)^p = \mu_p(y,F)^p = \lVert y\rVert_{S_p}^p,
\end{align*}
which proves that  $(e_{m_j})_{j\in\N}$ $1$-dominates $(u_j)_{j\in\N}$ in~$S_p$.

Before we tackle the case $E=B_p$, let us mention that Causey \cite[Lemma~3.2]{RC} has proved a closely related result, which states that $(e_{k_j})_{j\in\N}$ $4$-dominates $(u_j)_{j\in\N}$ in~$B_p$, where $k_j =  \min\operatorname{supp} u_j$ for $j\in\N$. Our result, with a larger constant, follows easily from Causey's. Nevertheless, we include a (fairly simple) proof of it to keep our work as self-contained as possible. 

Take $\mathcal{C}\in\operatorname{SC}$ such that $\lVert y\rVert_{B_p} = \beta_p(y,\mathcal{C})$. We may suppose that $F\subseteq\operatorname{supp} y$ for every $F\in\mathcal{C}$. Define $\mathcal{C}_j = \{ F\in\mathcal{C} : F\subseteq\operatorname{supp} u_j\}$ for $1\le j\le k$, $J = \{ j\in\{1,\ldots,k\} : \mathcal{C}_j\ne \emptyset\}$ and $\mathcal{C}' = \mathcal{C}\setminus\bigcup_{j\in J}\mathcal{C}_j$. Then~$\mathcal{C}$ is the disjoint union of $\{\mathcal{C}_j : j\in J\}\cup\{\mathcal{C}'\}$, so 
\begin{equation}\label{L:nbbs_Eq1}
    \lVert y\rVert_{B_p}^p = \sum_{j\in J}\beta_p(y,\mathcal{C}_j)^p + \beta_p(y,\mathcal{C}')^p.
\end{equation} 

It is is easy to estimate the first term on the right-hand side of this equation because the fact that $F\subseteq\operatorname{supp} u_j$ for each $F\in\mathcal{C}_j$ implies that
\[ \beta_p(y,\mathcal{C}_j) = \lvert\alpha_j\rvert\beta_p(u_j,\mathcal{C}_j)\le \lvert\alpha_j\rvert\,\lVert u_j\rVert_{B_p} = \lvert\alpha_j\rvert\qquad (j\in J). \]
Hence we have
\begin{equation}\label{L:nbbs_Eq2}
\sum_{j\in J}\beta_p(y,\mathcal{C}_j)^p\le \sum_{j\in J}\lvert\alpha_j\rvert^p\le \lVert x\rVert_{B_p}^p,\end{equation}
using the easy observation that the $B_p$-norm $1$-dominates the $\ell_p$-norm.

This completes the proof if $\mathcal{C}'=\emptyset$. Otherwise we can write $\mathcal{C}' = \{F_1<\cdots<F_n\}$ for some $n\in\N$. 
Set $K_r = \{ j\in\{1,\ldots,k\} : F_r\cap \operatorname{supp} u_j\ne\emptyset\}$ and $i_r = \min K_r$ for $1\le r\le n$. Then $M(K_r)\in\mathcal{S}_1$ because
\[ \lvert M(K_r)\rvert = \lvert K_r\rvert\le \lvert F_r\rvert\le\min F_r\le m_{i_r} = \min M(K_r).   \]
Furthermore, using that $\{F_r\cap \operatorname{supp} u_j : j\in K_r\}$ partitions~$F_r$, we obtain 
\begin{align}\label{L:nbbs_Eq4}
   \mu_1(y,F_r) &= \sum_{j\in K_r}\mu_1(y,F_r\cap\operatorname{supp}u_j) = \sum_{j\in K_r}\lvert\alpha_j\rvert\mu_1(u_j,F_r)\notag\\ &\le \sum_{j\in K_r}\lvert\alpha_j\rvert\,\lVert u_j\rVert_{B_p} = 
   \sum_{j\in K_r}\lvert\alpha_j\rvert = \mu_1(x,M(K_r)). 
\end{align}
The definition of~$\mathcal{C}'$ implies that $\lvert K_r\rvert\ge2$ and $K_r\setminus\{i_{r+1}\} < K_{r+1}$. In particular, we have $K_r<K_{r+2}$, so $M(K_r)<M(K_{r+2})$, and therefore we can define two Schreier chains by
\[ \mathcal{D} = \{ M(K_1)<M(K_3)<\cdots<M(K_{n'})\}\ \text{and}\ \mathcal{E} = \{M(K_2)<M(K_4)<\cdots<M(K_{n''})\}, \]
where $(n',n'') = (n-1,n)$ if~$n$ is even and $(n',n'') = (n,n-1)$ if~$n$ is odd. Hence, applying~\eqref{L:nbbs_Eq4} and using that $\mathcal{D}\cup\mathcal{E}$ partitions $\{ M(K_r) : 1\le r\le n\}$, we deduce that
\[ \beta_p(y,\mathcal{C}')^p = \sum_{r=1}^n\mu_1(y,F_r)^p\le \sum_{r=1}^n\mu_1(x,M(K_r))^p = \beta_p(x,\mathcal{D})^p + \beta_p(x,\mathcal{E})^p\le 2\lVert x\rVert_{B_p}^p. \]  
Finally, we substitute this upper bound together with that from~\eqref{L:nbbs_Eq2} into~\eqref{L:nbbs_Eq1} to conclude that $\lVert y\rVert_{B_p}^p\le 3\lVert x\rVert_{B_p}^p$, which proves that  $(e_{m_j})_{j\in\N}$ $3^{1/p}$-dominates $(u_j)_{j\in\N}$ in~$B_p$. 
 
\ref{L:nbbs2}. Set $N = \{ n_j:j\in\N\}$, and recall that $E_N = \overline{\operatorname{span}}\,\{e_{n_j} : j\in\N\}$. The $1$-un\-con\-di\-tion\-al\-i\-ty of the unit vector basis for~$E$ allows us to define a bounded operator $T\colon E\to E_N$ of norm at most~$\delta^{-1}$ by
\[ Tx = \sum_{k=1}^\infty\frac{\langle x, e_{n_k}^*\rangle}{\langle u_k,e_{n_k}^*\rangle}e_{n_k}. \]
It satisfies $Tu_j=e_{n_j}$ for each $j\in\N$ because $\langle u_j,e_{n_k}^*\rangle =0$ whenever $j\ne k$. Hence $(u_j)_{j\in\N}$ $\delta^{-1}$-dominates $(e_{n_j})_{j\in\N}$ in~$E$.

To prove that the subspace $W=\overline{\operatorname{span}}\{ u_j : j\in\N\}$ is $C_2$-complemented in~$E$, we note that $m_j = \max\operatorname{supp}u_j<n_{j+1}$ for each $j\in\N$ because $(u_j)_{j\in\N}$ is a block basic sequence. Therefore \Cref{C:intertwiningSeqs}\ref{C:intertwiningSeqs2} implies that the map $U\colon e_{n_j}\mapsto e_{m_j}$ for $j\in\N$ extends uniquely to an operator $U\in\mathscr{B}(E_N,E_M)$ with $\lVert U\rVert\le C_0$, where $C_0 = 2$ for $E=B_p$ and $C_0=2^{1/p}$ for $E=S_p$. By~\ref{L:nbbs1}, we can define an operator $V\in\mathscr{B}(E_M,W)$ with $\lVert V\rVert\le C_1$ by $Ve_{m_j} = u_j$ for each $j\in\N$. 
It follows that the composite operator $Q=VUT\in\mathscr{B}(E,W)$ satisfies $Qu_j = u_j$ for each $j\in\N$, so by linearity and continuity, $Q$ acts as the identity on~$W$; that is, $Q$ is a projection of~$E$ on\-to~$W$. Now the conclusion follows from the fact that
\[ \lVert Q\rVert \le\frac{C_0C_1}{\delta}\le \begin{cases} 
{\displaystyle{\frac{2\cdot 3^{1/p}}{\delta}}}\ &\text{for}\ E=B_p\\[7pt]
{\displaystyle{\frac{2^{1/p}}{\delta}}}\ &\text{for}\ E= S_p.  \end{cases} \]

\ref{L:nbbs3}, \ref{P:pointedbbs1}$\Rightarrow$\ref{P:pointedbbs2}. Suppose that $\inf_{j\in\N}\lVert u_j\rVert_\infty>0$, and define the sequences $(m_j)_{j\in\N}$ and $(n_j)_{j\in\N}$ as above. Then, as we already saw in the proof of~\ref{L:nbbs2}, the basic sequences $(u_j)_{j\in\N}$, $(e_{m_j})_{j\in\N}$ and $(e_{n_j})_{j\in\N}$ are all equivalent because $(u_j)_{j\in\N}$ dominates $(e_{n_j})_{j\in\N}$ by~\ref{L:nbbs2}, $(e_{n_j})_{j\in\N}$ dominates $(e_{m_j})_{j\in\N}$ by \Cref{C:intertwiningSeqs}\ref{C:intertwiningSeqs2}, and $(e_{m_j})_{j\in\N}$ dominates $(u_j)_{j\in\N}$ by~\ref{L:nbbs1}. 

The implication \ref{P:pointedbbs2}$\Rightarrow$\ref{P:pointedbbs3} is trivial.

Finally, we prove that~\ref{P:pointedbbs3} implies~\ref{P:pointedbbs1} by contradiction. Assume that $(u_j)_{j\in\N}$ dominates a subsequence of $(e_j)_{j\in\N}$ and that $\inf_{j\in\N}\lVert u_j\rVert_\infty=0$. \Cref{P:charDsubseq} shows that the latter assumption implies that $(u_j)_{j\in\N}$ admits a subsequence which is dominated by the unit vector basis $(d_j)_{j\in\N}$ for~$D$, where $D=c_0$ if $E=S_p$ and $D=\ell_p$ if $E=B_p$. Combining this with the former assumption, we conclude that $(d_j)_{j\in\N}$ dominates a subsequence of $(e_j)_{j\in\N}$, which is clearly impossible. 
\end{proof}

\section{The composition of two strictly singular operators is compact}\label{S:3}

\noindent The aim of this section is to prove \Cref{StimesSisK}. Our first step is to combine Propositions~\ref{P:charDsubseq}--\ref{L:nbbs} with the Principle of Small Perturbations to obtain the following result. 

\begin{lemma}\label{L:smallperturb}  
Let $(E,D) = (B_{p},\ell_{p})$ for some $1 < p < \infty$ or $(E,D) = (S_{p},c_{0})$ for some $1 \leq p < \infty$, and let $(w_n)_{n\in\N}$ be a weakly null sequence in~$E$ with $\inf_{n\in\N}\lVert w_n\rVert_E>0$.  
 \begin{enumerate}[label={\normalfont{(\roman*)}}]  
    \item\label{L:smallperturb1} Suppose that $\liminf_{n\to\infty}\lVert w_n\rVert_\infty=0$. Then  $(w_n)_{n\in\N}$ admits a subsequence $(w_n')_{n\in\N}$ which is a basic sequence equivalent to the unit vector basis for~$D$, and $\lVert w_n'\rVert_\infty\to 0$ as $n\to\infty$.
    \item\label{L:smallperturb2} Suppose that $\limsup_{n\to\infty}\lVert w_n\rVert_\infty>0$. Then $(w_n)_{n\in\N}$ admits a subsequence $(w_n')_{n\in\N}$ which is a basic sequence equivalent to a subsequence of the unit vector basis for~$E$, and $\inf_{n\in\N}\lVert w_n'\rVert_\infty>0$.
 \end{enumerate}
\end{lemma}

\begin{proof}
    We begin by noting that the weak convergence of $(w_n)_{n\in\N}$ implies that it is norm-bounded and hence semi\-normalized. For most of the proof, we shall consider the two cases simultaneously, but first we need to separate them:
    \begin{itemize}
        \item  In case~\ref{L:smallperturb1}, we replace $(w_n)_{n\in\N}$ with a subsequence such that $\lVert w_n\rVert_\infty\to0$ as $n\to\infty$, and we set $\xi = \inf_{n\in\N}\lVert w_n\rVert_E>0$.
        \item  In case~\ref{L:smallperturb2}, we replace $(w_n)_{n\in\N}$ with a subsequence such that $\xi := \inf_{n\in\N}\lVert w_n\rVert_\infty>0$ and observe that this implies that $\inf_{n\in\N}\lVert w_n\rVert_E\ge\xi$.
    \end{itemize}
    Returning to the unified approach, we define $m_0=0$, $P_0=0$ and $\varepsilon_j = \xi/(3\cdot 2^j+1)$ for each $j\in\N$. Using the fact that $(w_n)_{n\in\N}$ is weakly null, we can recursively construct in\-creasing sequences $(k_j)_{j\in\N}$ and $(m_j)_{j\in\N}$ of natural numbers such that the vectors $v_j := (P_{m_j}-P_{m_{j-1}})w_{k_j}\in E$ satisfy $\lVert w_{k_j} - v_j\rVert_E\le\varepsilon_j$ for each $j\in\N$, where~$P_m$ denotes the $m^{\text{th}}$ basis projection. 
    
    These choices imply that $(v_j)_{j\in\N}$ is a semi\-normalized block basic sequence because 
    \[  \sup_{n\in\N}\lVert w_n\rVert_E+\varepsilon_j\ge\lVert w_{k_j}\rVert_E+\lVert w_{k_j} - v_j\rVert_E\ge\lVert v_j\rVert_E \ge \lVert w_{k_j}\rVert_E-\lVert w_{k_j} - v_j\rVert_E\ge \xi-\varepsilon_j \] and $\varepsilon_j\le\xi/7$ for each $j\in\N$. Hence $(v_j)_{j\in\N}$ is equivalent to its normalization $(v_j/\lVert v_j\rVert_E)_{j\in\N}$.   
    Another application of the lower bound~$\xi-\varepsilon_j$ on~$\lVert v_j\rVert_E$ gives
    \[ \sum_{j=1}^\infty\frac{\lVert w_{k_j} - v_j\rVert_E}{\lVert v_j\rVert_E}\le \sum_{j=1}^\infty\frac{\varepsilon_j}{\xi-\varepsilon_j} = \sum_{j=1}^\infty\frac{1}{3\cdot 2^j} = \frac13<\frac12, \]
    so the Principle of Small Perturbations implies that $(w_{k_j})_{j\in\N}$ is a basic sequence equivalent to $(v_j)_{j\in\N}$ and therefore to $(v_j/\lVert v_j\rVert_E)_{j\in\N}$. We note that $\lVert w_{k_j} - v_j\rVert_\infty\le\varepsilon_j$ for each $j\in\N$.
    
    If we are in case~\ref{L:smallperturb1}, it follows that 
    \[ \frac{\lVert v_j\rVert_\infty}{\lVert v_j\rVert_E}\le\frac{\lVert w_{k_j}\rVert_\infty + \varepsilon_j}{6\xi/7}\to0\quad\text{as}\quad j\to\infty.  \]
    Combining this with \Cref{P:charDsubseq}, we deduce that $(v_j/\lVert v_j\rVert_E)_{j\in\N}$ admits a subsequence $(v_{j_i}/\lVert v_{j_i}\rVert_E)_{i\in\N}$ which is equivalent to the unit vector basis for~$D$, and therefore the same is true for the subsequence $(w_{k_{j_i}})_{i\in\N}$ of $(w_n)_{n\in\N}$. 

    Otherwise we are in case~\ref{L:smallperturb2}, and we see that $\lVert v_j\rVert_\infty\ge\lVert w_{k_j}\rVert_\infty - \varepsilon_j\ge 6\xi/7$, so 
    \[ \inf_{j\in\N}\frac{\lVert v_j\rVert_\infty}{\lVert v_j\rVert_E}\ge\frac{6\xi}{7\sup_{j\in\N}\lVert
     v_j\rVert_E}>0. \]
    Now \Cref{L:nbbs}\ref{L:nbbs3} implies that $(v_j/\lVert v_j\rVert_E)_{j\in\N}$ is equivalent to a subsequence of the unit vector basis for~$E$,  and therefore the same is true for the subsequence $(w_{k_j})_{j\in\N}$ of $(w_n)_{n\in\N}$.    
\end{proof} 

\begin{remark}
The two conditions in \Cref{L:smallperturb} are clearly not mutually exclusive, but at least one of them is always satisfied because the limit superior of a sequence dominates its limit inferior.
\end{remark}

\begin{lemma}\label{L:SSactiononWNS}
 Let~$T$ be a strictly singular operator on~$E$, where $E = B_{p}$ for some $1 < p < \infty$ or $E = S_{p}$ for some $1 \leq p < \infty$, and let $(w_n)_{n\in\N}$ be a weakly null sequence in~$E$ for which $\inf_{n\in\N}\lVert Tw_n\rVert_E>0$. Then $(w_n)_{n\in\N}$ admits a subsequence $(w_n')_{n\in\N}$ such that
\begin{enumerate}[label={\normalfont{(\roman*)}}]  
\item\label{L:SSactiononWNS:i} $\inf_{n\in\N} \lVert w_n'\rVert_\infty > 0$, and $(w_n')_{n\in\N}$ is equivalent to a subsequence of the unit vector ba\-sis for~$E$;
\item\label{L:SSactiononWNS:ii} $\lVert Tw_n'\rVert_\infty\to 0$ as $n\to\infty$, and $(Tw_n')_{n\in\N}$ is equivalent to the unit vector basis for~$D$, where $D=\ell_p$ if $E=B_p$ and $D=c_0$ if $E=S_p$.
\end{enumerate} 
\end{lemma}

\begin{proof} Let $(e_n)_{n\in\N}$ and $(d_n)_{n\in\N}$ denote the unit vector bases for~$E$ and~$D$, respectively. Since $\inf_{n\in\N}\lVert w_n\rVert_E\ge \lVert T\rVert^{-1}\inf_{n\in\N}\lVert Tw_n\rVert_E>0$, \Cref{L:smallperturb} implies that $(w_n)_{n\in\N}$ admits a subsequence $(w_n')_{n\in\N}$ which is a basic sequence and satisfies one of the following two conditions:  
\begin{enumerate}[label={\normalfont{(\Roman*)}}]  
\item\label{L:SSactiononWNS1} $(w_n')_{n\in\N}$ is equivalent to $(d_n)_{n\in\N}$, and $\lVert w_n'\rVert_\infty\to 0$ as $n\to\infty$; or 
\item\label{L:SSactiononWNS2} $(w_n')_{n\in\N}$ is equivalent to a subsequence of $(e_n)_{n\in\N}$, and $\inf_{n\in\N}\lVert w_n'\rVert_\infty>0$.
 \end{enumerate}
Being bounded, $T$ is weakly continuous, so the sequence $(Tw_n')_{n\in\N}$ converges weakly to~$0$, and therefore we may apply \Cref{L:smallperturb} to it, concluding that after replacing $(w_n')_{n\in\N}$ with a subsequence, we may suppose that $(Tw_n')_{n\in\N}$ is a basic sequence satisfying: 
\begin{enumerate}[label={\normalfont{(\Roman*)}},resume]
\item\label{L:SSactiononWNS3} $(Tw_n')_{n\in\N}$ is equivalent to $(d_n)_{n\in\N}$, and $\lVert Tw_n'\rVert_\infty\to 0$ as $n\to\infty$; or 
\item\label{L:SSactiononWNS4} $(Tw_n')_{n\in\N}$ is equivalent to a subsequence of $(e_n)_{n\in\N}$, and $\inf_{n\in\N}\lVert Tw_n'\rVert_\infty>0$.
 \end{enumerate}
Importantly, we note that conditions~\ref{L:SSactiononWNS1} and~\ref{L:SSactiononWNS2} remain true when we replace $(w_n')_{n\in\N}$ with a subsequence. 

We see that the desired conclusions~\ref{L:SSactiononWNS:i}--\ref{L:SSactiononWNS:ii} will follow if (and only if) conditions~\ref{L:SSactiononWNS2} and~\ref{L:SSactiononWNS3} are satisfied; we verify this by ruling out the other three possible combinations.

To present these arguments concisely, we require some notation. Let~$W$ and~$X$ denote the closed subspaces of~$E$ spanned by the basic sequences $(w_n')_{n\in\N}$ and $(Tw_n')_{n\in\N}$, respectively, and write~$\widetilde{T}$ for the restriction of~$T$ to~$W$, viewed as an operator into~$X$. We note that~$\widetilde{T}$ is strictly singular because~$T$ is. We introduce the following symbols for the isomorphisms that implement the equivalences in cases~\ref{L:SSactiononWNS1}--\ref{L:SSactiononWNS4}:
\begin{itemize}
\item In cases~\ref{L:SSactiononWNS1} and~\ref{L:SSactiononWNS3}, the linear maps~$U_1$ and $V_1$ given by $U_1d_j= w_j'$ and $V_1Tw_j'= d_j$ for $j\in\N$ define isomorphisms $U_1\in\mathscr{B}(D,W)$ and $V_1\in\mathscr{B}(X,D)$, respectively.   
\item In case~\ref{L:SSactiononWNS2}, we can take an infinite subset $M=\{m_1<m_2<\cdots\}$ of~$\N$ and an isomorphism $U_2\in\mathscr{B}(E_M,W)$ such that $U_2e_{m_j}=w_j'$ for each $j\in\N$, where we recall that $E_M = \overline{\operatorname{span}}\, \{e_{m_j} : j\in\N\}$. 
\item Similarly, in case~\ref{L:SSactiononWNS4}, we can take an infinite subset $N=\{n_1<n_2<\cdots\}$ of~$\N$ and an isomorphism $V_2\in\mathscr{B}(X,E_N)$ for which $V_2Tw_j'=e_{n_j}$ for each $j\in\N$. \end{itemize}

We are now ready to complete the argument by ruling out the three combinations of [\ref{L:SSactiononWNS1}~or~\ref{L:SSactiononWNS2}] and [\ref{L:SSactiononWNS3}~or~\ref{L:SSactiononWNS4}] that are not~\ref{L:SSactiononWNS2} and~\ref{L:SSactiononWNS3}: 
\begin{description}
  \item[\ref{L:SSactiononWNS1}+\ref{L:SSactiononWNS3}] In this case we have $V_1\widetilde{T}U_1d_j=d_j$ for each $j\in\N$, so $V_1\widetilde{T}U_1=I_D$, which contradicts that~$\widetilde{T}$ is strictly singular. 
  \item[\ref{L:SSactiononWNS1}+\ref{L:SSactiononWNS4}] In this case we see that $V_2\widetilde{T}U_1d_j = e_{n_j}$ for each $j\in\N$, which is impossible because $(d_j)_{j\in\N}$ does not dom\-i\-nate any subsequence of $(e_j)_{j\in\N}$.  
  \item[\ref{L:SSactiononWNS2}+\ref{L:SSactiononWNS4}]  In this case \Cref{C:intertwiningSeqs}\ref{C:intertwiningSeqs3} shows that the map $e_{m_{j_k}}\mapsto e_{n_{j_k}}$ for $k\in\N$ extends to an isomorphism of~$E_{M(J)}$ onto~$E_{N(J)}$ for some infinite subset $J=\{j_1<j_2<\cdots\}$ of~$\N$. However, this would imply that the restriction of~$V_2\widetilde{T}U_2$ to~$E_{M(J)}$ is an iso\-mor\-phism onto~$E_{N(J)}$, which again contradicts that~$\widetilde{T}$ is strictly singular. \qedhere
\end{description} 
\end{proof}

\begin{lemma}\label{L:easycompactnesscriterion}
    Let~$T$ be a non-compact operator from a Banach space~$X$ with a basis~$(x_n)_{n\in\N}$ into a Banach space~$Y$. Then~$(x_n)_{n\in\N}$ admits a normalized block basic sequence $(u_n)_{n\in\N}$ for which $\inf_{n\in\N}\lVert Tu_n\rVert_Y>0$.   
\end{lemma}

\begin{proof} Take $\eta\in(0,\lVert T+\mathscr{K}(X,Y)\rVert)$, where $\mathscr{K}(X,Y)$ denotes the closed subspace of $\mathscr{B}(X,Y)$ consisting of compact operators, and set $m_0=0$. We shall recursively choose natural numbers $m_1<m_2<\cdots$ and unit vectors $u_n\in\operatorname{span}\{ x_j : m_{n-1}<j\le m_n\}$ such that $\lVert Tu_n\rVert_Y\ge\eta/(K+1)$ for each $n\in\N$, where~$K$ denotes the basis constant of~$(x_n)_{n\in\N}$.

  We begin the construction by observing that since $\operatorname{span}\{ x_j : j\in\N\}$ is dense in~$X$ and $\lVert T\rVert>\eta$, we can find $m_1\in\N$ and a unit vector $u_1\in\operatorname{span}\{ x_j : 1\le j\le m_1\}$ such that $\lVert Tu_1\rVert_Y>\eta>\eta/(K+1)$. 

  Assume recursively that we have chosen natural numbers $m_1<m_2<\cdots<m_n$ and unit vectors $u_1,\ldots,u_n$ for some $n\in\N$. Since the basis projection~$P_{m_n}$ has finite rank, we have $\lVert T(I_X-P_{m_n})\rVert\ge \lVert T+\mathscr{K}(X,Y)\rVert>\eta$, so we can find $m_{n+1}>m_n$ and a unit vector $v_{n+1}\in\operatorname{span}\{ x_j : 1\le j\le m_{n+1}\}$ such that $\lVert T(I_X-P_{m_{n}})v_{n+1}\rVert_Y>\eta$. Then 
  \[ u_{n+1} = \frac{(I_X-P_{m_{n}})v_{n+1}}{\lVert (I_X-P_{m_{n}})v_{n+1}\rVert_X}\in\operatorname{span}\{ x_j : m_n< j\le m_{n+1}\} \]
  is a unit vector for which $\lVert Tu_{n+1}\rVert_Y\ge\eta/(K+1)$ because $\lVert (I_X-P_{m_{n}})v_{n+1}\rVert_X\le K+1$.  Hence the recursion continues, and the result follows. 
\end{proof}

\begin{proof}[Proof of Theorem~{\normalfont{\ref{StimesSisK}}}] 
  Let~$T$ and~$U$ be strictly singular operators on~$E$, where $E=B_p$ for some $1 < p < \infty$ or $E=S_{p}$ for some $1 \leq p < \infty$, and assume towards a contradiction that their composition~$TU$ is not compact. By \Cref{L:easycompactnesscriterion}, we can find a normalized block basic sequence $(u_n)_{n\in\N}$ of the unit vector basis $(e_n)_{n\in\N}$ for~$E$ such that $\inf_{n\in\N}\lVert TUu_n\rVert_E>0$. We note that $(u_n)_{n\in\N}$ is weakly null because $(e_n)_{n\in\N}$ is shrinking, so $(Uu_n)_{n\in\N}$ is also weakly null. Hence, applying \Cref{L:SSactiononWNS} with $w_n = Uu_n$ for $n\in\N$, we can extract a subsequence $(u_n')_{n\in\N}$ of $(u_n)_{n\in\N}$ for which 
  \begin{equation}\label{StimesSisK:eq1}
   \inf_{n\in\N} \lVert Uu_n'\rVert_\infty>0.    
  \end{equation}
  Since $(u_n')_{n\in\N}$ is weakly null and $\inf_{n\in\N} \lVert Uu_n'\rVert_E\ge\inf_{n\in\N} \lVert Uu_n'\rVert_\infty>0$, another application of \Cref{L:SSactiononWNS} shows that $(u_n')_{n\in\N}$ admits a subsequence $(u_n'')_{n\in\N}$ such that $\lVert Uu_n''\rVert_\infty\to0$ as $n\to\infty$. However, this contradicts~\eqref{StimesSisK:eq1}.
\end{proof}

Having presented a self-contained proof of \Cref{StimesSisK}, we shall next outline how it may alternatively be deduced from the work of Androulakis, Dodos, Sirotkin and Troit\-sky~\cite{ADST} concerning compositions of strictly singular operators. This argument still relies on Lemmas~\ref{L:smallperturb} and~\ref{L:SSactiononWNS}, so it is not independent of the above proof. It is based on the observation that every strictly singular operator on the Baernstein and Schreier spaces is ``finitely strictly singular'' in the following sense, a result that may be of interest in its own right. 

\begin{definition}\label{D:FSS} An operator $T\in\mathscr{B}(X,Y)$ between Banach spaces~$X$ and~$Y$ is \emph{finitely strictly singular} if, for every $\varepsilon>0$, there exists $n\in\N$ such that every subspace of~$X$ of dimension at least~$n$ contains a unit vector~$x$ for which $\lVert Tx\rVert\le\varepsilon$.    
\end{definition}

\begin{proposition}\label{charSS}
    Let $(E,D) = (B_{p},\ell_{p})$ for some $1 < p < \infty$ or $(E,D) = (S_{p},c_{0})$ for some $1 \leq p < \infty$. The following conditions are equivalent for an operator $T\in\mathscr{B}(E)\colon$
\begin{enumerate}[label={\normalfont{(\alph*)}}] 
\item\label{charSS1} $T$ is finitely strictly singular;
\item\label{charSS2} $T$ is strictly singular;
\item\label{charSS3} $T$ does not fix a copy of~$D;$
\item\label{charSS4} the identity operator on~$D$ does not factor through~$T$. 
\end{enumerate}
\end{proposition}    

The implications \ref{charSS1}$\Rightarrow$\ref{charSS2}$\Rightarrow$\ref{charSS3}$\Rightarrow$\ref{charSS4} are clear, and~\ref{charSS4} implies~\ref{charSS2} because~$E$ is saturated with complemented copies of~$D$ by \cite[Theorem~2.4]{LS}, so it only remains to prove that~\ref{charSS2} implies~\ref{charSS1}. This relies on Milman's ``flat-vector lemma'', originally published in~\cite{milman}; we refer to \cite[Lemma~13]{CP} for an easily accessible presentation in English. 

\begin{lemma}[Milman]\label{L:milman} Let $W$ be a non-zero subspace of~$c_0$, and take a natural number $n\le\dim W$. Then~$W$ contains a unit vector~$w$ which attains its norm in at least~$n$ coordinates; that is, $\lVert w\rVert_\infty=1$ and the set $\{ j\in\N : \lvert w(j)\rvert=1\}$ has cardinality at least~$n$.     
\end{lemma}

\begin{proof}[Proof of Proposition~{\normalfont{\ref{charSS}}}, \ref{charSS2}$\Rightarrow$\ref{charSS1}] Assume towards a contradiction that $T\in\mathscr{B}(E)$ is a strictly singular operator which fails to be finitely strictly singular. Then, for some $\varepsilon>0$, there is a sequence $(W_n)_{n\in\N}$ of subspaces of~$E$ such that
\begin{equation}\label{charSS:eq1}
   \dim W_n\ge 2n-1\qquad\text{and}\qquad \lVert Tw\rVert_E\ge\varepsilon\lVert w\rVert_E\qquad (n\in\N,\,w\in W_n). 
\end{equation}
As previously observed, $E$ is a vector subspace of~$c_0$ because the unit vector basis is a normalized basis for~$E$. Hence \Cref{L:milman} applies, showing that for each $n\in\N$, we can choose $v_n\in W_n$ such that $\lVert v_n\rVert_\infty=1$ and the set $J_n = \{ j\in\N : \lvert v_n(j)\rvert=1\}$ has cardinality at least~$2n-1$.  Then $J_n\cap[n,\infty)$ has cardinality at least~$n$, so~$J_n$ contains a Schreier set~$F_n$ of cardinality~$n$, and therefore $\lVert v_n\rVert_{S_p}\ge \mu_p(v_n,F_n)=n^{1/p}$ for every $p\in[1,\infty)$. Since $\lVert\,\cdot\,\rVert_{B_p}\ge\lVert\,\cdot\,\rVert_{S_1}$, we conclude that $\lVert v_n\rVert_E\to\infty$ as $n\to\infty$. Consequently, defining $w_n = v_n/\lVert v_n\rVert_E\in W_n$ for each $n\in\N$, we have 
\begin{equation}\label{charSS:eq2} 
\lVert w_n\rVert_\infty = 1/\lVert v_n\rVert_E\to 0\quad\text{as}\quad n\to\infty. 
\end{equation}
In particular, $(w_n)_{n\in\N}$ is a norm-bounded sequence with $\langle w_n, e_j^*\rangle\to0$ as $n\to\infty$ for each $j\in\N$, so $(w_n)_{n\in\N}$ is weakly null because the basis $(e_j)_{j\in\N}$ for~$E$ is shrinking. Since $\inf_{n\in\N}\lVert Tw_n\rVert_E\ge\varepsilon$ by~\eqref{charSS:eq1}, \Cref{L:SSactiononWNS} implies that $(w_n)_{n\in\N}$ admits a subsequence $(w_n')_{n\in\N}$ for which $\inf_{n\in\N}\lVert w_n'\rVert_\infty>0$. However, this contradicts~\eqref{charSS:eq2}.
\end{proof}

\begin{remark}
  In the case $E=S_p$, $1\le p<\infty$, conditions \ref{charSS1}--\ref{charSS3} in \Cref{charSS} remain equivalent for operators $T\in\mathscr{B}(S_p,Y)$ into any Banach space~$Y$. The implications \ref{charSS1}$\Rightarrow$\ref{charSS2}$\Rightarrow$\ref{charSS3} are true in general, and~\ref{charSS3} implies~\ref{charSS2} because~$S_p$ is saturated with copies of~$c_0$, as before. To see that~\ref{charSS2} implies~\ref{charSS1}, we proceed as in the above proof, noting that we do not use the facts that $T$ is strictly singular or that the codomain of~$T$ is~$E$ until we invoke \Cref{L:SSactiononWNS} in the penultimate line. 
  For an operator $T\in\mathscr{B}(E,Y)$, still assumed not to be finitely strictly singular, but now with arbitrary co\-domain~$Y$, we can instead apply \Cref{L:smallperturb}\ref{L:smallperturb1}: since $(w_n)_{n\in\N}$ is weakly null with $\lVert w_n\rVert_{E}=1$ and $\lVert w_n\rVert_\infty\to0$ as \mbox{$n\to\infty$}, it admits a subsequence $(w_n')_{n\in\N}$ which is equivalent to the unit vector basis $(d_n)_{n\in\N}$ for~$D$. Taking $U\in\mathscr{B}(D,E)$ with $Ud_n=w_n'$ for each $n\in\N$, we have 
  \[ \inf_{n\in\N}\lVert TUd_n\rVert_Y=\inf_{n\in\N}\lVert Tw_n'\rVert_Y\ge\varepsilon. \] If $D=c_0$, this condition implies that the restriction of~$TU$ to the closed subspace of~$c_0$ spanned by \mbox{$\{ d_n : n\in N\}$} is an isomorphic embedding for some infinite subset~$N$ of~$\N$ by a famous result of Rosen\-thal, originally stated as the first remark following \cite[Theorem~3.4]{ROS}; that is, $T$ fixes a copy of~$c_0$. 
\end{remark}

To present our alternative proof of \Cref{StimesSisK} based on the results of~\cite{ADST}, we require the following definitions from that paper; throughout, $X$~denotes a Banach space.
    \begin{itemize}
    \item A semi\-normalized basic sequence $(x_n)_{n\in\N}$ in~$X$ is \emph{Schreier spreading} if there is a constant $C\ge 1$ such that, for every pair of non-empty Schreier sets $F,G\in\mathcal{S}_1$ of the same cardinality~$n$, the finite basic sequences $(x_j)_{j\in F}$ and $(x_j)_{j\in G}$ are $C$-equivalent; that is, writing $F=\{ f_1<f_2<\cdots<f_n\}$ and $G = \{ g_1<g_2<\cdots<g_n\}$, we have
    \begin{equation}\label{Defn:schreierspreading} \frac{1}{C}\biggl\|\sum_{j=1}^n\alpha_j x_{f_j}\biggr\|_X\le\biggl\|\sum_{j=1}^n\alpha_j x_{g_j}\biggr\|_X\le C\biggl\|\sum_{j=1}^n\alpha_j x_{f_j}\biggr\|_X\qquad (\alpha_1,\ldots,\alpha_n\in\K). \end{equation}
    The set of semi\-normalized  basic sequences in~$X$ that are Schreier spreading and weakly null is denoted $\operatorname{SP}_{1,w}(X)$.
    \item Given two semi\-normalized Schreier spreading basic sequences $(x_n)_{n\in\N}$ and $(y_n)_{n\in\N}$ in~$X$, define $(x_n)_{n\in\N}\approx_1(y_n)_{n\in\N}$ if there is a constant $C\ge 1$ such that the finite basic sequences $(x_n)_{n\in F}$ and $(y_n)_{n\in F}$ are $C$-equivalent for every $F\in\mathcal{S}_1$.     
    This defines an equivalence relation~$\approx_1$ on the set of semi\-normalized Schreier spreading basic sequences in~$X$.
\end{itemize}
We shall repeatedly use the following simple facts about these notions. 
\begin{itemize} 
\item Equivalence of basic sequences implies $\approx_1$-equivalence; that is, suppose that $(x_n)_{n\in\N}$ and $(y_n)_{n\in\N}$ are semi\-normalized Schreier spreading basic sequences in~$X$ which are equivalent. Then $(x_n)_{n\in\N}\approx_1(y_n)_{n\in\N}$. 
\item Suppose that $(x_n')_{n\in\N}$ is a subsequence of a semi\-normalized Schreier spreading basic sequence $(x_n)_{n\in\N}$ in~$X$. Then $(x_n')_{n\in\N}$ is Schreier spreading, and $(x_n')_{n\in\N}\approx_1 (x_n)_{n\in\N}$; see \cite[Proposition~3.5(i)]{ADST}.
\end{itemize}

\begin{lemma}\label{L:approx1rel}
Let $E = B_{p}$ for some $1 < p < \infty$ or $E = S_{p}$ for some $1 \leq p < \infty$.
\begin{enumerate}[label={\normalfont{(\roman*)}}] 
\item\label{L:approx1rel:i} The unit vector basis for~$E$ belongs to~$\operatorname{SP}_{1,w}(E)$.
\item\label{L:approx1rel:ii} Suppose that $(x_n)_{n\in\N},(y_n)_{n\in\N}\in\operatorname{SP}_{1,w}(E)$ satisfy either 
\begin{equation}\label{L:approx1rel:eq1}  \liminf_{n\to\infty}\lVert x_n\rVert_\infty=\liminf_{n\to\infty}\lVert y_n\rVert_\infty=0 
\end{equation}
or 
\begin{equation}\label{L:approx1rel:eq2}  \min\Bigl\{\limsup_{n\to\infty}\lVert x_n\rVert_\infty, \limsup_{n\to\infty}\lVert y_n\rVert_\infty\Bigr\}>0. 
\end{equation}
Then $(x_n)_{n\in\N}\approx_1(y_n)_{n\in\N}$. 
\end{enumerate}
\end{lemma}

\begin{proof} \ref{L:approx1rel:i}. We already know that the unit vector basis~$(e_n)_{n\in\N}$ for~$E$ is normalized and shrinking, and thus weakly null, so it remains only to verify~\eqref{Defn:schreierspreading}. However, for $F\in\mathcal{S}_1$ and scalars $\alpha_j\in\K$, $j\in F$, we have
\[ \biggl\lVert\sum_{j\in F}\alpha_j e_j\biggr\rVert_E = \begin{cases} {\displaystyle{\sum_{j\in F} \lvert \alpha_j\rvert}}\ &\text{if}\ E=B_p\\ 
{\displaystyle{\Bigl(\sum_{j\in F} \lvert \alpha_j\rvert^p\Bigr)^{1/p}}}\ &\text{if}\ E=S_p,
\end{cases} \]
which shows that~\eqref{Defn:schreierspreading} is satisfied for $C=1$.

\ref{L:approx1rel:ii}. \Cref{L:smallperturb} applies because the sequences $(x_n)_{n\in\N}$ and $(y_n)_{n\in\N}$ are semi\-normalized and weakly null. Hence, if~\eqref{L:approx1rel:eq1} is satisfied, $(x_n)_{n\in\N}$ and $(y_n)_{n\in\N}$ admit subsequences~$(x_n')_{n\in\N}$ and $(y_n')_{n\in\N}$ which are equivalent to the unit vector basis for~$D$, where $D=\ell_p$ if $E=B_p$ and $D= c_0$ if $E=S_p$. It follows that $(x_n')_{n\in\N}$ and $(y_n')_{n\in\N}$ are equivalent to each other, and using the two bullet points above, we obtain
    \[ (x_n)_{n\in\N}\approx_1 (x_n')_{n\in\N}\approx_1 (y_n')_{n\in\N}\approx_1 (y_n)_{n\in\N}.  \]

    On the other hand, if~\eqref{L:approx1rel:eq2} is satisfied, $(x_n)_{n\in\N}$ and $(y_n)_{n\in\N}$ admit subsequences $(x_n')_{n\in\N}$ and $(y_n')_{n\in\N}$ which are equivalent to subsequences $(e_{j_n})_{n\in\N}$ and $(e_{k_n})_{n\in\N}$ of the unit vector basis $(e_n)_{n\in\N}$ for~$E$. Since $(e_n)_{n\in\N}\in\operatorname{SP}_{1,w}(E)$ by~\ref{L:approx1rel:i}, the two bullet points above imply that
    \[ (x_n)_{n\in\N}\approx_1 (x_n')_{n\in\N}\approx_1 (e_{j_n})_{n\in\N}\approx_1(e_n)_{n\in\N}\approx_1 (e_{k_n})_{n\in\N}\approx_1 (y_n')_{n\in\N}\approx_1 (y_n)_{n\in\N}. \]
    
    In both cases, the conclusion that $(x_n)_{n\in\N}\approx_1(y_n)_{n\in\N}$ follows from transitivity of~$\approx_1$. 
\end{proof}

\begin{proof}[Alternative proof of Theorem~{\normalfont{\ref{StimesSisK}}}, based on Androulakis et al~\cite{ADST}] We check that the conditions for applying the ``Moreover'' statement in the first part of \cite[Theorem~4.1]{ADST} are satisfied for $\xi=1$ and $n=2$: 
\begin{itemize}
    \item \Cref{L:approx1rel} implies that $\operatorname{SP}_{1,w}(E)$ contains at most two distinct $\approx_1$-equivalence classes because every sequence $(x_j)_{j\in\N}$ in~$E$ trivially satisfies \[ \limsup_{j\to\infty}\lVert x_j\rVert_\infty\ge\liminf_{j\to\infty}\lVert x_j\rVert_\infty\ge 0. \] 
    \item Combining \Cref{charSS} with \cite[Proposition~2.4(i)]{ADST}, we see that every strictly singular operator on~$E$ is ``$\mathcal{S}_1$-strictly singular'' in the terminology of \cite[Definition~2.1]{ADST}. 
    \item No subspace of~$E$ is isomorphic to~$\ell_1$ because~$E$ has a shrinking basis.
\end{itemize}
Hence \cite[Theorem~4.1]{ADST} shows that the composition of two strictly singular operators on~$E$ is compact.
\end{proof}

\section{Subsymmetric basic sequences in the Baernstein and Schreier spaces}\label{S:4}
\noindent A basis for a Banach space is \emph{subsymmetric} if it is unconditional and equivalent to all its subsequences. The unit vector bases for~$c_{0}$ and~$\ell_{p}$, for $1 \leq p < \infty$, are standard examples of subsymmetric bases. Using our previous results, we can show that every subsymmetric basic sequence in the Baernstein and Schreier spaces is equivalent to one of those. 

This relies on the following lemma, which is certainly known to specialists; we include a short, elementary  proof for ease of reference. 
\begin{lemma}\label{L:bonus}
    Let $(x_n)_{n\in\N}$ be a norm-bounded unconditional basis for a Banach space~$X$. Then either $(x_n)_{n\in\N}$ is weakly null or it admits a subsequence which is equivalent to the unit vector basis for~$\ell_1$. 
\end{lemma}

\begin{proof}
  Suppose that $(x_n)_{n\in\N}$ is not weakly null.  Then there is a functional $x^*\in\ X^*$ of norm~$1$ such that the sequence $(\langle x_n,x^*\rangle)_{n\in\N}$ does not tend to~$0$, so we can pass to a subsequence $(x_{n_j})_{j\in\N}$ for which $\eta := \inf_{j\in\N}\lvert \langle x_{n_j},x^*\rangle\rvert >0$. 

  Given $m\in\N$ and $\alpha_1,\ldots,\alpha_m\in\K$, for each $1\le j\le m$, we can choose $\sigma_j,\tau_j\in\K$ with $\lvert\sigma_j\rvert = \lvert\tau_j\rvert = 1$ such that $\sigma_j \langle x_{n_j},x^*\rangle\ge\eta$ and $\tau_j\alpha_j\ge 0$. Let $K$ denote the unconditional constant of the basis $(x_n)_{n\in\N}$. Then we have
  \begin{align*}
      K\biggl\lVert\sum_{j=1}^m\alpha_j x_{n_j}\biggr\rVert &\ge \biggl\lVert\sum_{j=1}^m\sigma_j\tau_j\alpha_j x_{n_j}\biggr\rVert\ge \biggl\lvert\Bigl\langle \sum_{j=1}^m\sigma_j\tau_j\alpha_j x_{n_j}, x^*\Bigr\rangle\biggr\rvert\\ &= \biggl\lvert \sum_{j=1}^m\tau_j\alpha_j \sigma_j\langle x_{n_j}, x^*\rangle\biggr\rvert = \sum_{j=1}^m\lvert\alpha_j\rvert\,\lvert \langle x_{n_j}, x^*\rangle\rvert\ge \eta \sum_{j=1}^m\lvert\alpha_j\rvert.
  \end{align*}
  This shows that the subsequence $(x_{n_j})_{j\in\N}$ $K/\eta$-dominates the unit vector basis for~$\ell_1$, which on the other hand trivially dominates every bounded sequence by the triangle inequality. The result follows.  
  \end{proof}

\begin{remark}\label{R:subsymm}
    A subsymmetric basis $(x_n)_{n\in\N}$ must be semi\-normalized. Indeed, if it had no lower norm bound, we could recursively choose an increasing sequence $1\le m_1<m_2<\cdots$ of integers such that $\lVert x_{m_n}\rVert\le \lVert x_{n}\rVert/n$ for each $n\in\N$. However, this would contradict that $(x_{m_n})_{n\in\N}$ dominates $(x_{n})_{n\in\N}$. A similar argument shows that $\sup_{n\in\N}\lVert x_n\rVert<\infty$. 
\end{remark}

\begin{proposition}\label{characterisesubsymmetric}
    Let $(E,D) = (B_{p},\ell_{p})$ for some $1 < p < \infty$ or $(E,D) = (S_{p},c_{0})$ for some $1 \leq p < \infty$. A basic sequence in~$E$ is subsymmetric if and only if it is equivalent to the unit vector basis for~$D$.
\end{proposition}

\begin{proof}
    The backward implication is clear. For the forward implication, suppose that $(u_{n})_{n\in\N}$ is a subsymmetric basic sequence in~$E$. \Cref{R:subsymm} shows that it is semi\-nor\-mal\-ized, so combining \Cref{L:bonus} with the fact that no subspace of~$E$ is isomorphic to~$\ell_1$ because~$E$ has a shrinking basis, we conclude that $(u_{n})_{n\in\N}$ is weakly null. Hence, by \Cref{L:smallperturb}, it admits a subsequence $(u_{k_{n}})_{n\in\N}$ which is equivalent to either the unit vector basis for~$D$ or a subsequence $(e_{m_{n}})_{n\in\N}$ of the unit vector basis for~$E$. In the former case, the result follows because $(u_{n})_{n\in\N}$ is equivalent to $(u_{k_{n}})_{n\in\N}$. The latter case is impossible because it would imply that $(e_{m_{n}})_{n\in\N}$ is subsymmetric, contradicting \cite[Theorem~4.1]{LS}. 
\end{proof}

\section{Uncomplemented block subspaces of the Baernstein and Schreier spaces}\label{S:5}

\noindent The aim of this section is to show that the Baernstein and Schreier spaces contain block basic sequences whose closed span is not complemented. We follow the approach that Gasparis and Leung took when proving \cite[Proposition~4.7]{GL}, which is the counterpart for the higher-order Schreier spaces~$X[\mathcal{S}_n]$ for $n\in\N$. It is based on a lemma of Lindenstrauss and Tzafriri and involves the following standard notion. 
\begin{definition}
    Let~$X$ be a Banach space with a basis~$(x_n)_{n\in\N}$. A block basic sequence $(u_n)_{n\in\N}$ of $(x_n)_{n\in\N}$ is \emph{skipped} if there are integers $0=m_0<m_1<m_2<\cdots$ such that \[ u_n\in\operatorname{span}\{ x_j : m_{n-1}<j<m_n\}\qquad (n\in\N). \]  
\end{definition}

\begin{proposition}\label{P:uncomplementedblock}
The Baernstein spaces~$B_p$, for $1<p<\infty$, and the Schreier spaces~$S_p$, for $1\le p<\infty$, contain block basic sequences whose closed span is not complemented. 

More precisely, let $(E,D) = (B_p,\ell_p)$ for some $1<p<\infty$ or $(E,D) = (S_p,c_0)$ for some $1\le p<\infty$, and let $(u_n)_{n\in\N}$ be a skipped, normalized block basic sequence of the unit vector basis for~$E$ such that the unit vector basis for~$D$ dominates $(u_n)_{n\in\N}$. For each $n\in\N$, take $m_n\in(\max\operatorname{supp} u_n, \min\operatorname{supp} u_{n+1})\cap\N$, and set
\[ t_n = \frac1k\qquad (n\in[2^{k-1},2^k)\cap\N,\, k\in\N). \]
Then $(u_n + t_ne_{m_n})_{n\in\N}$ is a block basic sequence whose closed span is not complemented in~$E$. 
\end{proposition}

\begin{proof}
    We can choose~$m_n$ as specified because the block basic sequence $(u_n)_{n\in\N}$ is skipped, and it ensures that $(u_n + t_ne_{m_n})_{n\in\N}$ is a block basic sequence. Assume towards a contradiction that its closed span is complemented in~$E$. Since $t_n\to 0$ as $n\to\infty$, \cite[Lemma~2.a.11]{LT1} implies that the basic sequence $(u_n)_{n\in\N}$  dominates $(t_ne_{m_n})_{n\in\N}$. By hy\-poth\-e\-sis, the unit vector basis $(d_n)_{n\in\N}$ for~$D$ dominates $(u_n)_{n\in\N}$, and $(t_ne_{m_n})_{n\in\N}$ dominates $(t_ne_n)_{n\in\N}$ by \Cref{C:intertwiningSeqs}\ref{C:intertwiningSeqs1}. Hence we can find a constant $C>0$ such that $(d_n)_{n\in\N}$ $C$-dominates $(t_ne_n)_{n\in\N}$, so in particular, we have
    \begin{equation}\label{P:uncomplementedblock:eq}
        C\biggl\lVert\sum_{n=2^{k-1}}^{2^k-1}d_n\biggr\rVert_D\ge \biggl\lVert\sum_{n=2^{k-1}}^{2^k-1}t_ne_n\biggr\rVert_E = \frac1k\biggl\lVert\sum_{n=2^{k-1}}^{2^k-1}e_n\biggr\rVert_E
    \end{equation}
    for every $k\in\N$. However, 
    \[ \biggl\lVert\sum_{n=2^{k-1}}^{2^k-1}d_n\biggr\rVert_D = \begin{cases} 1\ &\text{for}\ D=c_0\\ 2^{\frac{k-1}{p}}\ &\text{for}\ D=\ell_p \end{cases}\qquad\text{and}\qquad \biggl\lVert\sum_{n=2^{k-1}}^{2^k-1}e_n\biggr\rVert_E = \begin{cases} 2^{\frac{k-1}{p}}\ &\text{for}\ E=S_p\\ 2^{k-1}\ &\text{for}\ E=B_p  \end{cases} \]
    because $[2^{k-1},2^k)\cap\N\in\mathcal{S}_1$. Substituting these values into~\eqref{P:uncomplementedblock:eq}, we conclude that
    \[ C\ge \begin{cases} {\displaystyle{\frac{2^{\frac{k-1}{p}}}{k}}}\ &\text{for}\ (E,D)=(S_p,c_0)\\[1.5ex]
    {\displaystyle{\frac{\bigl(2^{1-\frac{1}{p}}\bigr)^{k-1}}{k}}}\ &\text{for}\ (E,D)=(B_p,\ell_p), \end{cases} \] 
    which is absurd because the right-hand sides are unbounded as $k\to\infty$ in both cases.
\end{proof}

\begin{remark} 
\begin{enumerate}[label={\normalfont{(\roman*)}}] 
\item In order to apply \Cref{P:uncomplementedblock}, we must find a skipped, normalized block basic sequence $(u_n)_{n\in\N}$ of the unit vector basis for~$E$ that is dominated by the unit vector basis for~$D$. \Cref{P:charDsubseq} ensures that many such sequences exist.  
\item One reason that \Cref{P:uncomplementedblock} is significant is that in the case of the Baernstein spaces, it disproves \cite[Lemma~II.3.1]{S}, in which Seifert stated that for every \mbox{$1< p<\infty$}, every closed subspace of~$B_p$ spanned by a semi\-normalized block basic sequence is complemented in~$B_p$. (We note that ``semi\-normalized'' is irrelevant here because a block basic sequence $(u_n)_{n\in\N}$ and its normalization $(u_n/\lVert u_n\rVert)_{n\in\N}$ span the same subspace.) 
\end{enumerate}
\end{remark}

We conclude with a dichotomy which recombines some of our previous results to show that although the closed span of a block basic sequence need not be complemented in the Baernstein or Schreier spaces, their block basic sequences are nevertheless intimately connected to complemented subspaces.
 
\begin{proposition}\label{P:complementationDichotomy} Let $(E,D) = (B_p,\ell_p)$ for some $1<p<\infty$ or $(E,D) = (S_p,c_0)$ for some $1\le p<\infty$, and let $(u_n)_{n\in\N}$ be a normalized block basic sequence of the unit vector basis for~$E$. 
\begin{enumerate}[label={\normalfont{(\roman*)}}]
  \item\label{P:complementationDichotomy:i} If $\inf_{n\in\N}\lVert u_n\rVert_\infty>0$, then $(u_n)_{n\in\N}$ is equivalent to a subsequence of the unit vector basis for~$E$, and $\overline{\operatorname{span}}\, \{u_n : n\in\N\}$ is complemented in~$E$.
  \item\label{P:complementationDichotomy:ii} Otherwise $(u_n)_{n\in\N}$ admits a subsequence $(u_{n_j})_{j\in\N}$ that is equivalent to the unit vector basis for~$D$ and $\overline{\operatorname{span}}\, \{u_{n_j} : j\in\N\}$ is complemented in~$E$.
\end{enumerate}
\end{proposition}

\begin{proof}
    Part~\ref{P:complementationDichotomy:i} follows from \Cref{L:nbbs}\ref{L:nbbs2}--\ref{L:nbbs3}. 

    To prove~\ref{P:complementationDichotomy:ii}, suppose that $\inf_{n\in\N}\lVert u_n\rVert_\infty=0$. As in the proof of \cite[Proposition~2.14]{LS}, we see that $(u_n)_{n\in\N}$ admits a subsequence $(u_{n_j})_{j\in\N}$ that is dominated by the unit vector basis $(d_j)_{j\in\N}$ for~$D$ by \cite[Lemma~2.8]{LS}, so we have an operator $U\in\mathscr{B}(D,E)$ such that $Ud_j = u_{n_j}$ for $j\in\N$. On the other hand, \cite[Lemma~2.10]{LS} shows that there is an operator $V\in\mathscr{B}(E,D)$ such that $Vu_{n_j}=d_j$ for $j\in\N$. Hence $(u_{n_j})_{j\in\N}$ is equivalent to $(d_j)_{j\in\N}$, and~$UV$ is a projection of~$E$ onto $\overline{\operatorname{span}}\, \{u_{n_j} : j\in\N\}$.
\end{proof}

\begin{acknowledgements} This work is part of the second-named author's PhD at Lan\-caster Uni\-ver\-sity. He acknowledges with thanks the funding from the EPSRC (grant reference code EP/V520214/1) that has supported his studies.
\end{acknowledgements}

\end{document}